\documentclass[preprint,11pt,authoryear]{article}

\usepackage{amssymb, amstext,latexsym, amsthm, amsmath, amscd, amsfonts}
\usepackage[english]{babel}
\usepackage[latin1]{inputenc}
\usepackage[all]{xy}
\usepackage{natbib}

\usepackage[normalem]{ ulem }
\usepackage{soul}

\usepackage{epsfig, float}
\usepackage{graphicx}
\usepackage{mathtools}
\usepackage{color}
\usepackage{pgf,pgfarrows,pgfnodes}
\usepackage{xmpmulti}
\usepackage{caption}

\usepackage{tikz}
\usetikzlibrary{calc,arrows,decorations.pathmorphing,intersections}
\usetikzlibrary{patterns}

\definecolor{vert_fonce}{rgb}{0,0.75,0.25}
\definecolor{vert_pire}{rgb}{0,0.70,0.30}
\definecolor{vert_plus}{rgb}{0,0.60,0.40}
\definecolor{rouge_plus}{rgb}{0.70,0.20,0.10}
\definecolor{bleu_plus}{rgb}{0.30,0,0.70}

\RequirePackage[colorlinks=true,citecolor=blue,linkcolor=blue]{hyperref}
\usepackage{hyperref}

\newtheorem{ex}{Example}[section]
\newtheorem{theorem}{Theorem}[section]
\newtheorem{prop}[theorem]{Proposition}

\newtheorem{rem}[theorem]{Remark}

\newcommand{\beq} {\begin{eqnarray*}}
\newcommand{\eeq} {\end{eqnarray*}}
\newcommand{\trm} {\textrm}
\newcommand{\tbf} {\textbf}

\def\QED{\hfill\vrule height 1.5ex width 1.4ex depth -.1ex \vskip20pt}

\def\1{\mathbf{1}}

\def \R{\mathbb R}

\def \Z{{\mathbf Z}}
\def \Y{{\mathbf Y}}
\def \E{{\mathbb E}}
\def \Var{\hbox{{\rm Var}}}

\def \Cov{\hbox{{\rm Cov}}}

\def \leqp{\leqslant}
\def \geqp{\geqslant}

\usepackage{bbm}					

\newcommand{\ind}{\mathbbm{1}}

\newcommand{\diff}{\mathop{}\mathopen{}\mathrm{d}}

\begin{document}


\title{Multilevel branching splitting algorithm for estimating rare event probabilities}

\maketitle

\begin{center}
Agn\`{e}s Lagnoux and Pascal Lezaud
\end{center}

\textbf{Abstract:} \small{We analyse the splitting algorithm performance in the estimation of rare event probabilities in a discrete multidimensional framework. For this we assume that each threshold is partitioned into disjoint subsets and the probability for a particle to reach the next threshold will depend on the starting subset. A straightforward estimator of the rare event probability is given by the proportion of simulated particles for which the rare event occurs. The variance of this estimator we get is the sum of two parts: one part resuming the variability due to each threshold and a second part resuming the variability due to the thresholds number. This decomposition is analogous to that of the continuous case. The optimal algorithm is then derived by cancelling the first term leading to optimal thresholds. 
Then we compare this variance with that of the algorithm in which one of the threshold has been deleted. 
Finally, we investigate the sensitivity of the variance of the estimator with respect to a shape deformation of an optimal threshold. As an example, we consider a two-dimensional Ornstein-Uhlenbeck process with conformal maps for shape deformation.}

\textbf{Keywords:} \small{splitting; rare event probability estimation;  Monte Carlo;  branching process;  simulation;  variance reduction;  first crossing time density;  conformal maps}

\section{Introduction}

The risk modelling approach consists in firstly formalizing the system
considered and secondly using mathematical or simulation tools to obtain some
estimates \citep{Aldous89,Sadowsky96}. Analytical and numerical approaches are useful, but may
require many simplifying assumptions. On the other hand, Monte Carlo simulation
is a practical alternative when the analysis calls for fewer simplifying
assumptions. Nevertheless, obtaining accurate estimates of rare event
probabilities, say about $10^{-9}$ to  $10^{-12}$, using traditional techniques require a huge amount
of computing time.

Many techniques for reducing the number of trials in Monte Carlo simulation have
been proposed, like importance sampling or
trajectory splitting (\citet{LLL09}). In the splitting technique, we suppose there exists some well identifiable intermediate 
states that are visited
much more often than the target states themselves and behave as gateways to
reach the rare event. Thus we consider a
decreasing sequence of events $B_i$ leading to the rare event $B$:
\begin{align}\label{def:emb}
B\coloneqq B_{M+1}\subset B_M \subset \ldots \subset B_1\,.
\end{align}
Then $p \coloneqq \mathbb{P}(B)=\mathbb{P}(B|B_M)\mathbb{P}(B_M|B_{M-1})\ldots \mathbb{P}(B_2|B_1)\mathbb{P}(B_1)$ where on the right hand side, each conditioning event is "not rare". These conditional probabilities are in general not available
explicitly. Instead, we know how to make evolve the particles from level $B_i$
to the next level $B_{i+1}$ (e.g. Markovian behaviour). 

The principle of the algorithm is at first to run  
simultaneously several particles starting from the level $B_i$;
after a while, some of them have evolved "badly", the other have evolved "well" i.e. have succeeded in reaching the  
threshold $B_{i+1}$. Then "bad"
particles are moved to the position of the "good" ones and so on until $B$ is
reached. In such a way, the more promising 
particles are favoured. Examples of this class of algorithms can be found in \citet{Aldous-Vazirani94} with the "go with
the winners" scheme, in \citet{Jerrum-Sinclair97} and \citet{Diaconis-Holmes95} in 
approximate counting and in a more general setting in \citet{Doucet-Freitas-Gordon01, DelMoral04, Cerou-Guyader05, DMG05, MBJV14} . 

The difficulty comes from the complexity of the dynamics of the particles. A simpler analysis can be done focusing only on the 
underlying Markov chain that represents the changes of thresholds. In this technique, we make a Bernoulli trial to check whether or not 
the set event $B_1$ has occurred. In that case, we split this trial in $R_1$ Bernoulli
subtrials and for
each of them we check again whether or not the event $B_2$ has occurred. This procedure is repeated at each
level, until $B$ is reached. If an event level is not reached, neither is $B$, then we stop the current
retrial. Using $N$ independent replications of this procedure, we have then
considered $N R_1 \ldots R_M$ trials, taking into account for example, that if we have failed to reach a
level $B_i$ at the $i$-th step, the $R_i \ldots R_M$ possible retrials have
failed. Clearly the particles reproduce and evolve independently. 

An unbiased estimator of $p$ is given by the quantity
\begin{equation*}
\widehat p_{M+1}=\frac{N_B}{N \prod_{i=1}^M R_i},
\end{equation*}
where $N_B$ is the total number of trajectories having reached the set
$B$. Considering that this algorithm is represented by $N$ independent Galton-Watson branching
processes, as done in \citet{Lagnoux06}, the variance of $\widehat p_{M+1}$ can then be derived and 
depends on the probability transitions and on the mean numbers of particles successes at each level. 
Leading by the heuristic presented in \citet{Villen91b,Villen97}, an optimal algorithm is derived by minimising the variance 
of the estimator for a given budget (or computational
cost). This cost is defined as the expected number of trials generated during the simulation,
each trial being weighted by a cost function.

The optimisation of the algorithm suggests to take all the transition probabilities 
equal to a constant and the numbers of splitting equal to the inverse of this constant \citet{Lagnoux06}. Then we deduce the number of thresholds $M$ and finally the number $N$ of replications. In fact, optimal values are chosen in 
such a way to balance between the increase of the variance when the number splitting is small and the
exponential growth in computational effort when too much splitting are used.

In higher dimension, the engineering community have proposed algorithms  to estimate rare event probabilities. Subset simulation which is also based on a partitioning of the space into nested subsets uses Markov Chain simulation (in particular the Metropolis Hastings scheme) \cite{AB01}.  Importance sampling techniques have also been developed in that framework. When the failure region is not too complex to describe, schemes to construct importance sampling algorithms have been introduced that are based on design points  (see e.g. \cite{APB99,DKD98} and the references therein) or adaptive pre-samples (see e.g. \cite{AB99} and the references therein). When the complexity of the rare event increases, it seems to be difficult to construct efficient importance sampling scheme \cite{SPP93}.

In this paper, we continue the multidimensional approach and study theoretically the algorithm introduced in \citet{GHSZ98} and \citet{Garvels00} mainly in order to obtain a new expression of the variance of the estimator ana\-lo\-gous to that of the continuous case (\citet{LLL09}). Thus, we assume that each threshold is partitioned into $s$ disjoint subsets and the probability for a particle starting from a threshold to reach the next threshold will depend on the starting subset. Unlike the unidimensional case, the hardness to reach the next threshold differs according to the starting subset; in some sense the threshold is no longer an iso-probability level. In this context, the variance of the estimator $\widehat p_{M+1}$ is the sum of two parts: one part resuming the variability due to each threshold and a second part resuming the variability due to the thresholds number (see Proposition \ref{prop:var}). For the unidimensional case, only the second term remains. The optimal algorithm is then derived by cancelling the first term of the variance leading to iso-probability levels and by optimising the other parameters as in the unidimensional case. 

Furthermore, by introducing new operators, we obtain an alternative expression of the variance which is more tractable when we wish to compare the variance of the estimators in an algorithm with $M$ thresholds with the variance in an algorithm in which one of the threshold has been deleted. 
More precisely, we study the need of an intermediate threshold and derive a procedure to detect whether we shall keep it or not. In order to obtain a simple criteria, we assume the optimal shape of the thresholds of the optimal algorithm. Finally, we investigate the sensitivity of the variance of $\widehat{p}_{M+1}$ with respect to a shape deformation of the threshold relatively to the optimal shape. 

The remainder of this paper is divided into five sections. In Sections 2--4 we present, analyse theoretically and optimise the splitting algorithm in the multidimensional case. Then, Sections~\ref{sec:form_it_var} and \ref{sec:perturb} deal with the sensitivity analysis of the variance as previously presented. In particular, in Section \ref{sec:perturb}, we illustrate a way to deform the shape of the thresholds to get uniform occupation densities with a 2D Ornstein-Uhlenbeck process. Finally, we complete the paper by a conclusion and some perspectives. More details and all the proofs are postponed in the appendices.

\section{Multilevel Splitting Algorithm}

\subsection{Definition of the thresholds and related tools}

In order to estimate the probability $p$ that a particle starting from a point in some state space $E$ reaches the critical subset $B\subset E$, we use the so-called splitting algorithm based on the nested sequence $B_1,\ldots,B_{M+1}$ defined in~\eqref{def:emb}.
Moreover, each frontier $\partial B_k$ of $B_k$ is partitioned into $s$ disjoint subsets, denoted $\partial B_k^{(i)}$, such that
\[\partial B_k=\bigcup_{i=1}^s \partial B_k^{(i)}, \quad k=1,\ldots, M.\]

We assume that each $\partial B_k$ has the same number $s$ of subsets; this assumption is not restrictive as one can see in the sequel. In any case, one can obviously rewrite the problem
under concern in this particular setting.

The random dynamics of the particle are modelled by a stochastic process $Y=(Y_t; t\geqp 0)$ and for $k=1,\ldots,M+1$, we define $\tau_k$ as the first time that the particle hits $\partial B_k$. Hence $p$ can be written as $p=\mathbb{P}(\tau_{M+1}<\infty)$. For the sake of simplicity, we assume naturally that $Y$ evolves continuously and all the intermediate thresholds are hit if the last one is. In fact, the dynamics under concern is not directly the particle one but rather the one of the embedded Markov chain observed at each time the particle hits a frontier $\partial B_k$. This embedded Markov chain will be denoted $(X_k)_{0\leqp k \leqp M+1}$. Thus, $X_k=i$ if the particle at time $\tau_k$ lies in $\partial B_k^{(i)}$ i.e. $Y_{\tau_k} \in\partial B_k^{(i)}$.

\paragraph{Measures $\gamma_k$ and functions $f_k$} We define for any $k=1,\ldots, M$, a measure $\gamma_k$ on the frontier $\partial B_k$ by
\[\gamma_k(i)=\mathbb{P}(X_k=i\;;\; \tau_k<\infty).\]
This measure acts on the functions $f$ defined on $\partial B_k$ by $\gamma_k(f)=\mathbb{E}\left[f(X_k)\;;\; \tau_k<\infty\right]$
in such a way that $\gamma_k(\mathbf{1})=\mathbb{P}(\tau_k<\infty)$ is the probability that the particle hits the event $B_k$ ($\mathbf{1}$ stands for the unit function). 

For any $k=1,\ldots, M$, we denote $\mathcal M_k$ (resp. $\mathcal F_k$) the set of measures (resp. functions) defined on $\partial B_k$. In particular, the functions $f_k \in \mathcal F_k$ defined by
\[f_k(i)=\mathbb{P}(\tau_{M+1} <\infty \mid X_k=i \; ; \;  \tau_k <\infty), \quad k=1,\ldots, M\]
play a special role, since 
\begin{equation}\label{eq:p}
\gamma_k(f_k)=\sum_{i=1}^s \gamma_k(i)f_k(i)=p, \quad k=1,\ldots, M.
\end{equation}
In fact, $f_k(i)$ quantifies the hardness to reach the target set $B$ starting from $\partial B_k^{(i)}$ while $\gamma_k(i)f_k(i)$ quantifies the hardness to reach $B$ passing by $\partial B_k^{(i)}$ and starting from $O$. Furthermore, for $k=2,\ldots, M$, we introduce the operators $P_k$, $k=2,\cdots,M$ defined on $\partial B_{k-1} \times \partial B_{k}$ by
$P_k(i,j)=\mathbb{P}(X_k =j\; ; \; \tau_k<\infty| X_{k-1}=i\; ; \; \tau_{k-1} <\infty)$. 
Nevertheless it is easier to consider $P_k$ as an operator right acting on $\mathcal F_k$ as an operation $\mathcal F_k\to \mathcal F_{k-1}$ according to
\[
P_k(f)(i)=\mathbb{E}\left[f(X_k)\;;\; \tau_k<\infty \mid X_{k-1}=i \; ; \; \tau_{k-1} <\infty \right]
\]
and left acting on $\mathcal M_{k-1}$ as an operation $\mathcal M_{k-1}\to \mathcal M_{k}$ according to $(\mu P_k)(f)=\mu(P_k f)$.

Each operator $P_k$ is not Markovian, since the probability to reach $\partial B_k$ is not equal one; hence we define $g_{k-1}\in \mathcal F_{k-1}$ by
\begin{equation}\label{def:G}
g_{k-1}(i)\coloneqq P_k(\mathbf{1})(i)=\mathbb{P}(\tau_k<\infty\mid X_{k-1}=i\; ; \; \tau_{k-1} <\infty),
\end{equation}
for $k=2, \ldots, M$. Remark that there is no need to define $g_M$ since it would correspond to $f_M$.

We easily get the following transport relations for $k=2,\ldots, M$,
\begin{equation}\label{eq:dyngamma}
\gamma_k=\gamma_{k-1} P_k,\quad f_{k-1}= P_k(f_k).
\end{equation}

\begin{figure}[h!]
\centering
\includegraphics[height=10cm]{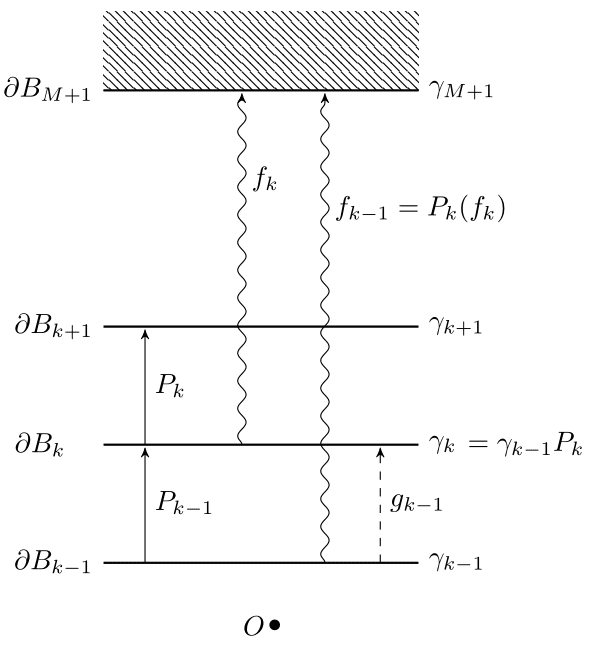}
\caption{This figure summarizes the notation previously introduced.}\label{fig:notation}
\end{figure}

The notation is summarized in Figure~\ref{fig:notation}.

\noindent\\
\paragraph{Normalized measures $\mu_k$} Since $\gamma_k$ is not a probability measure, we define its normalized version $\mu_k$ on $\partial B_k$ that acts on the functions $f \in \mathcal F_k$ in the following way
\begin{equation*}
\mu_k(f)\coloneqq \frac{\gamma_k(f)}{\gamma_k(\mathbf{1})}=\mathbb{E}\left[f(X_k) \mid\tau_k <\infty \right],
\end{equation*}
(assuming that the thresholds have been chosen such that $\gamma_k(\mathbf{1})\neq 0$ for all $k$). 
We notice that 
\begin{equation}\label{eq:muk_gk}
\mu_k(g_k)=\frac{\gamma_{k+1}(\mathbf{1})}{\gamma_k(\mathbf{1})}=\mathbb{P}(\tau_{k+1}<\infty\mid\tau_k<\infty)
\end{equation}
and
\begin{equation*}
\mu_k(f_k)=\frac{p}{\gamma_k(\mathbf{1})}=\mathbb{P}(\tau_{M+1}<\infty\mid\tau_k<\infty).
\end{equation*}

Equation~\eqref{eq:dyngamma} induces the following scheme for the dynamics of  $\mu_k$
\begin{equation} \label{eq:dynmu}
\mu_k=\frac{\gamma_{k-1}(\mathbf{1})}{\gamma_k(\mathbf{1})} \mu_{k-1}P_k=\frac{1}{\mu_{k-1}(g_{k-1})} \mu_{k-1}P_k
\end{equation}
that leads to
\begin{equation*}
\gamma_k(f)=\mu_k(f) \gamma_1(\mathbf{1})\prod_{p=1}^{k-1}\mu_p(g_p),
\end{equation*}
which applied successively to the functions $f_k$ and $f_{k+1}$ yields to $\mu_k(f_k)=\mu_{k+1}(f_{k+1}) \mu_k(g_k).$


\paragraph{Convention} We extend the previous definitions to $k=0$ and $k=M+1$. Considering that the particles are generated at the same point $O$, we define $\mathcal F_0$ as the set of constant functions and in particular $f_0=p$ and $g_0=\gamma_1(\mathbf{1})$. Analogously $\mathcal M_0$ will represent the set of the Dirac measures at $O$ up to a constant. Hence $\gamma_0(f)=f$ (and $\mu_0=\gamma_0$). Obviously, $\gamma_0(\mathbf{1})=\mathbb P(\tau_0<\infty)=1$.  

In the same way, $B_{M+1}$ is reduced to a unique point, denoted e.g. by $\omega$. Then $\mathcal F_{M+1}$ is reduced to the constant functions, with $f_{M+1}=1$ and $\mathcal M_{M+1}$ is the set of the Dirac measures at $\omega$ up to a constant, with $\gamma_{M+1}(f)=f p$ (such as $\gamma_{M+1}(f_{M+1})=\gamma_{M+1}(\mathbf{1})=p$) and $\mu_{M+1}(f)=f$. 
We set also $P_1(f)=\gamma_1(f)$ and  $P_{M+1}(f)=f \times f_M$.

\subsection{Multilevel Splitting Algorithm}\label{sec:algo}

To estimate the rare event probability we proceed according to the algorithm already introduced in \citet{GHSZ98} and \citet{Garvels00}. Its principle is the following:

\begin{description}
\item[Initialization:] We perform independently $N$ particles from the same starting point $O$. A random number $Z_1$ of particles reach the threshold $B_1$, where $Z_1$ has a binomial distribution with parameters $N$ and $\gamma_1(\mathbf{1})$. These $Z_1$ particles are spread over the subsets $\partial B_1^{(i)}$ according to a multinomial random variable (r.v.) $\text{Mult}(Z_1, \mu_1)$. Let $\Z_1$ be the corresponding random vector $(Z_{1 1},\ldots, Z_{1 r})$.

\item[Step $\bf n$ $(2 \leqp n \leqp  M)$:] Each of the $Z_{n-1}$ particles in $\partial B_{n-1}$ is duplicated $R_{n-1}$ times; so that a total number  $R_{n-1} Z_{n-1}$ of particles is achieved. These new particles evolve according to the dynamics of the original process and the number $Z_{n j}$ of particles reaching $\partial B_n^{(j)}$ is still a random number. Consider now the random vector $\Z_n=(Z_{n 1}, \ldots, Z_{n r})$. The $Z_{n j}$ particles in $B_n^{(j)}$ come from different subsets $\partial B_{n-1}^{(i)}$; then we decompose $Z_{n j}$ in the following sum 
\begin{equation}\label{eq:somme}
Z_{n j}=\sum_{i=1}^s Y_{nj}^{i},
\end{equation}
where $Y_{nj}^{i}$ is the number of particles from $\partial B_{n-1}^{(i)}$ and having reached $\partial B_n^{(j)}$ whose total number 
$Y_{n}^i=\sum_{j=1}^s Y_{n j}^{i}$ is a binomial r.v. with parameters $R_{n-1} Z_{(n-1)i}$ and $g_{n-1}(i)$.

We represent the numbers $Y_{nj}^i$ in a $s \times s$ tabular where each line $\Y_n^i=(Y_{n 1}^i,\ldots, Y_{n s}^i)$, conditionally to the knowledge of the total number $Y_n^i$, is distributed as a multinomial r.v. with parameters $Y_n^i$ and  $Q_n(i,\cdot)$ where 
\[Q_n(i,\cdot)\coloneqq \frac{P_n(i,\cdot)}{g_{n-1}(i)}=\mathbb{P}(X_n=\cdot \mid X_{n-1}=i \;;\; \tau_n<\infty).\]
\begin{table}[hbt]
\begin{center}
\begin{tabular}{|c|c|c||c|}
\hline
$\cdot$&$\cdot$&$\cdot$&$\cdot$\\
\hline
$\cdot$& $Y_{nj}^i$&$\cdot$&$Y_n^i$\\
\hline
$\cdot$&$\cdot$&$\cdot$&$\cdot$\\
\hline
\hline
$\cdot$&$Z_{nj}$&$\cdot$&$Z_n$\\
\hline
\end{tabular}
\end{center}
\end{table}
In a nutshell, the random vector $\Z_n$ can be expressed as the sum $\Z_n=\sum_{i=1}^s \Y_n^i$ of the random vectors $\Y_n^i$ and the total number of particles at the end of step $n$ is $Z_n=\sum_{i,j} Y_{nj}^i=\sum_{i=1}^s Y_n^i$.

\item[Final step:] Each of the $Z_M$ particles in $\partial B_M$ is duplicated $R_M$ times to get a total number $R_M Z_M$ of particles.  These new particles evolve accordingly to the dynamics of the original process and the $Z_{M+1}$ particles having reached $\partial B_{M+1}$ come from different subsets $\partial B_{M}^{(i)}$; then we decompose $Z_{M+1}$ in the following sum 
\begin{equation}\label{eq:2_somme}
Z_{M+1}=\sum_{i=1}^s Y_{M+1}^{i}
\end{equation}
where $Y_{M+1}^{i}$ represents the number of particles from $\partial B_{M}^{(i)}$ and having reached $\partial B_{M+1}$. Conditionally to the random vector $\Z_M$, the r.v.s $Y_{M+1}^i$, $i=1,\ldots, r$ are independent and distributed as a binomial r.v. with parameters $R_M Z_{M i}$ and $ f_M(i)$. The set $B_{M+1}$ being reduced to a point, the result of this final step is simply the total number $Z_{M+1}$ of particles in $B_{M+1}$.
\end{description}

\section{Algorithm analysis}\label{sec:algo_an}

In this section, we present a natural unbiased estimator of $p$ and give several expressions of its variance including the one given in \cite{GHSZ96}. We also define the cost of the algorithm.

\subsection{A natural unbiased estimator of $p$}\label{sec:est}

An estimator of the probability to hit $\partial B_{n+1}$ conditionally that $\partial B_n$ has been hit is naturally given by the ratio between the number of particles in $\partial B_{n+1}$ and $R_n$ times the number of particles in $\partial B_n$, from which we deduce a natural estimator of the probability of interest $p$ 
\begin{equation}\label{eq:est}
\widehat{p}_{M+1}=\frac{Z_{1}}{N} \times \prod_{n=1}^{M-1} \frac{Z_{n+1}}{R_n Z_{n}} \times \frac{Z_{M+1}}{R_M Z_{M}}=\frac{Z_{M+1}}{N R_1\ldots R_M}.
\end{equation}

Introducing the deterministic quantities $r_0=N$ and $r_n=R_n r_{n-1}$, $n=1,\ldots,M$, leads to 
$\widehat{p}_{M+1}=Z_{M+1}/r_M$.
%
Then it is obvious to show that this estimator is unbiased. Indeed by conditioning,~\eqref{eq:2_somme} yields 
\[\E[\widehat{p}_{M+1}]=\frac{1}{r_{M-1}}\sum_{i=1}^s  \E[Z_{M i}]f_M(i).\]
To derive the mean of $Z_{M i}$, notice that $\E[Z_{1j}]=N\gamma_1(j)$, $j=1,\ldots,r$. By a new conditioning, we get that 
\[\E[Z_{2j}]=R_1 \sum_{i=1}^s \E[Z_{1i}] P_2(i,j)=N R_1 \gamma_2(j)=r_1 \gamma_2(j),\]
the last equality coming from~\eqref{eq:dyngamma}. An induction principle allows us to establish that for $n=2,\ldots, M$,
\begin{equation}\label{eq:moypar}
\E[Z_{nj}]=r_{n-1} \gamma_n(j)
\end{equation}
that leads to $\E[\widehat{p}_{M+1}]= \sum_{j=1}^s \gamma_M(j) f_M(j)=\gamma_M(f_M)=p$.

\subsection{The variance of the estimator}

\begin{prop}\label{prop:var} The coefficient of variation is given by
\begin{equation}\label{eq:var_mes}
\frac{\Var(\widehat{p}_{M+1})}{p^2}=\sum_{k=1}^M \frac{1}{\gamma_k(\mathbf{1})}\left(\frac{1}{r_{k-1}}-\frac{1}{r_{k}}\right)\frac{\Var_{\mu_k}(f_k)}{\E^2_{\mu_k}(f_k)} + \sum_{k=0}^M \frac{1}{r_k  \gamma_k(\mathbf{1})}\frac{1 - \mu_k(g_k)}{\mu_k(g_k)}.
\end{equation}

Introducing the operators $\Gamma_{i+1}$ defined, for $f$, $g\in \mathcal{F}_{i+1}$, by
$\Gamma_{i+1}(f,g)=P_{i+1}(fg)-P_{i+1}(f) P_{i+1}(g)$, we have
\begin{equation}\label{eq:var14}
\Var(\widehat{p}_{M+1})=\sum_{i=0}^{M} \frac{1}{r_i } \gamma_i\left(\Gamma_{i+1}(f_{i+1})\right).
\end{equation}
\end{prop}


The variance is then split into two parts. The first sum outlines the variability due to the shape of the thresholds $\partial B_k$ (defined by the $f_k$'s) whereas the second outlines the variability due to the thresholds number $M$, replication numbers $R_k$ and thresholds position (contained in the $P_k$'s and $g_k$'s).
Also we refer to \ref{app:op_Gamma} for more details on the operators $\Gamma_{i+1}$. 

\paragraph{Comparison with other algorithms} Notice that for $s=1$, the measures $\gamma_k$ and the functions $f_k$ are constant.
Since $\mu_k(g_k)=\gamma_{k+1}/\gamma_k$, the expression of the variance becomes
\[\frac{\Var(\widehat{p}_{M+1})}{p^2}=\sum_{k=0}^M \frac{1 - \mu_k(g_k)}{r_k \gamma_{k+1}}=\sum_{k=0}^M
\frac{1}{r_{k}}\left(\frac{1}{\gamma_{k+1}}-\frac{1}{\gamma_{k}}\right)\]
that corresponds to the expression established in \citet{Lagnoux08}.

Furthermore, Formula~\eqref{eq:var_mes} corresponds to equation (2.21) established in \citet{Garvels00} for an algorithm with a single intermediate threshold. It also has been established in \citet{LLL09,CDMG2011} in the general and continuous settings. 

Finally, simple computation leads to the following expression
\begin{equation*}
\frac{\Var(\widehat{p}_{M+1})}{p^2}=
\sum_{k=1}^{M} \frac{1}{\gamma_k(\mathbf{1})}\left(\frac{1}{r_{k-1}}-\frac{1}{r_k}\right) \frac{\mu_k(f_k^2)}{\mu_k^2(f_k)} + \left(\frac{1}{pr_M}-\frac{1}{r_0}\right)
\end{equation*}
that can be found in \citet{GHSZ99}.

\subsection{The cost of the algorithm}

The efficiency of the algorithm can be traduced in terms of the variance of the estimator that must be the smallest possible under the condition that the cost (in terms of computer time for example) remains finite. Our goal is then to derive the optimal parameters of the algorithm for a fixed cost.

The total number of particles generated during the algorithm is the r.v.
$N+ R_1 Z_1 + \ldots +R_M Z_M$. From~\eqref{eq:moypar}, $\E[Z_n]=r_{n-1} \gamma_n(\mathbf{1})$,
the mean of the total number of particles generated by the algorithm is  
\begin{equation*}
C_{M+1}^{(0)} \coloneqq r_0 + r_1 \gamma_1(\mathbf{1}) + \ldots +r_M \gamma_M(\mathbf{1})
\end{equation*}
and can be considered as a natural cost.

Now we present a more realistic cost that takes into account the probability $P_k(i,j)$ to reach $\partial B_k^{(j)}$ from $\partial B_k^{(i)}$. Actually, even if the algorithm presented here is based on the simulation of multinomial r.v.s, the introduction of this new cost allows to consider the dynamics of a particle between two successive thresholds through the functions $g_k$. Thus we associate to each particle from $\partial B_{k}^{(i)}$ a unitary cost $c_{k}(i)$ that depends on the starting threshold and the  hardness $g_{k}(i)$ to succeed in reaching the next threshold. More precisely, we assume that
\[c_0=c(\gamma_1(\mathbf{1})),  \quad c_{k}(i)=c(g_{k}(i)), \quad k=1,\ldots, M,\]
where $c$ is a positive function, decreasing (the smaller the probability of success is the highest the cost is) such that $c(x)$ converges to a constant (in general small) when $x$ tends to $1$.

\begin{prop}\label{prop:cost}
The mean cost is given by 
\begin{equation}\label{eq:cout}
C_{M+1}=N c_0+\sum_{n=1}^{M} r_{n} \sum_{i=1}^s  \gamma_{n}(i)c_{n}(i)=\sum_{n=0}^{M}r_{n} \gamma_{n}(c_{n}).
\end{equation}
\end{prop}

The approach presented here leads to a relatively simple formula for the total mean cost, similar to the one used in \citet{Lagnoux08}. The multidimensionality of the model is taken into account through the function $c_n$.

\section{Algorithm optimisation}\label{sec:algo_opt}
 
Before proceeding to the optimisation of the algorithm, we start recalling the general setting.

\subsection{General setting}\label{ssec:gen_setting}

In many applications, the rare event probability $p$ can be viewed as an overflow probability. More precisely, let $h$ be a real-valued measurable function defined on $E$ and $L\geqp 0$ be a given threshold. Then $p$ is rewritten as $p=\mathbb{P}(h(Y_t)\geqp L)$ where the process $Y$ has been defined in the Introduction. As a consequence, we can naturally use the function $h$ to determine the intermediate thresholds and apply the splitting methodology to the real-valued process $Z$ defined by $Z_t\coloneqq h(Y_t)$, for all $t\geqp 0$ (for simplicity $Z_0\geqp 0$). For the sake of simplicity we assume that $Z$ evolves continuously and all the intermediate thresholds are hit if the last one is. 

However remark that the intermediate thresholds $L_1$, \ldots, $L_{M+1}$ for $Z$ define splitting surfaces $\partial B_1$, \ldots, $\partial B_{M+1}$ for $Y$ by $\partial B_k=\{y\in E \mid h(y)=L_k\}$. Defining the levels $\partial B_k$ in such a way is not well adapted and is far to be optimal. Indeed, this methodology is geometrical and only based on a level set without taking into account the probabilistic aspects. More precisely, it seems natural to incorporate information of the hardness to reach the target set from any point of the $\partial B_k$. This information is precisely given by the function $f_k$ introduced previously. So, assuming the possibility to define a function $f$, named \emph{importance function}, on the whole space $E$ by 
\[
f(x)=\mathbb{P}(\tau_{M+1}<\infty \mid \text{starting from } x),
\] 
we rather define $\partial B_k$ as the set of the points $x\in E$ such that $f(x)=L_k$ for some $L_k \in [0,1]$. In some sense, we use iso-probability density levels as intermediate thresholds. Of course, the difficulty here is to determine the function $f$, the thresholds number $M$ and the values of $L_k$. Nevertheless, there exists methods that allow to get estimators of $f$ using a reverse time analysis as proposed in \citet{Garvels00}. 

To illustrate the importance of a good choice of the intermediate thresholds, let us consider the following example represented in Figure \ref{fig:dim2}.

\begin{ex}
With $M=1$ and a threshold $\partial B_1$ partitioned in two subsets such that
\[\left\{\begin{array}{ll}
&\gamma_1(1)= 10^{-2},\;  \gamma_1(2)=0.5, \\
&f_1(\mathbf{1})=10^{-1},\;  f_1(2)=10^{-3},\\
\end{array}\right.\]
we obtain $p=1.5 \cdot{} 10^{-3}$ and $\gamma_1(\tbf{1})=0.51$.

Let us simulate particles starting from $O$. We expect that $51\%$ of them reach the threshold $\partial B_1$, with $50\%$ in $\partial B_1^{(2)}$ and only $1\%$ in $\partial B_1^{(1)}$. Nevertheless among those in $\partial B_1$, the particles in the subset $\partial B_1^{(1)}$ have $100$ times more likely to reach the target set than those in $\partial B_1^{(2)}$. So, using this design of $\partial B_1$ leads to simulate almost $50\%$ of particles pointlessly. We see that using a function $f_1=p/\alpha$, with $\alpha \in ]0,1[$, implies that $\gamma(\tbf{1})=\alpha$ and all the particles in $\partial B_1$ then have the same probability to reach the target set. 

\begin{figure}[h!]
\centering
\includegraphics[height=6cm]{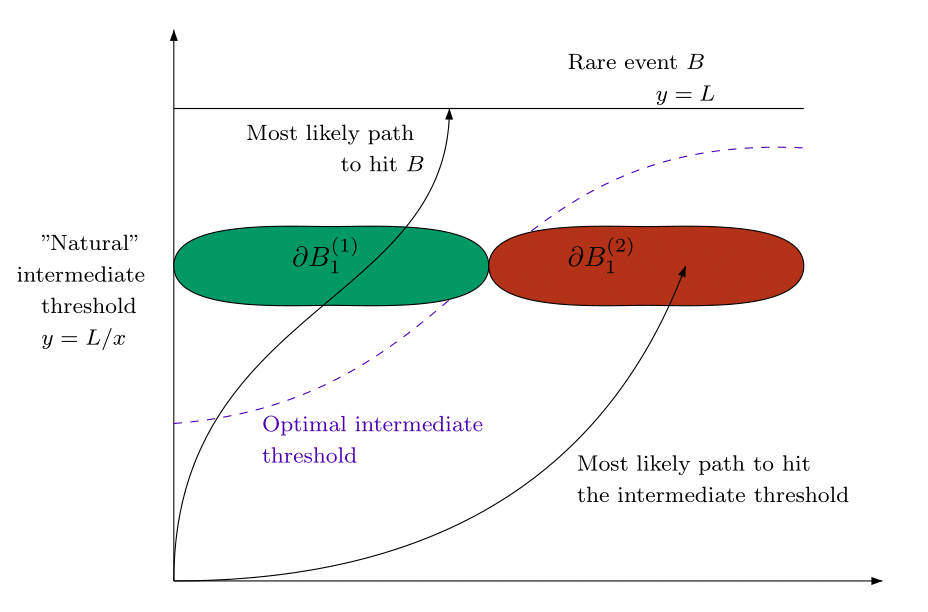}
\caption{The crucial choice of the importance function on an example.}\label{fig:dim2}
\end{figure}
\end{ex}

The construction of the importance function when the target probability has a large deviation characterization is handled in \citet{DD07}. This context is also considered in \citet{Sadowsky96} and Remark~\ref{rm:Mir}. Nevertheless, it seems difficult to translate the results obtained into the framework of this paper.

\subsection{Optimisation}\label{ssec:opt}

It is important to keep in mind Equation~\eqref{eq:var_mes} and the fact that the variance of the $\widehat{p}_{M+1}$ can be split in two parts: a first one 
resuming the variability due to the shape of the thresholds and a second one
resuming the variability due to the thresholds number, replication numbers and thresholds position.

Furthermore, the splitting algorithm's parameters are: the initial number $N$ of particles, the replication numbers $R_1,\ldots,R_M$, the number $M$ of intermediate thresholds and their characteristics (through the $P_k$'s and the $g_k$'s). 

\begin{prop}\label{prop:opt}
The parameters of the algorithm optimised by minimisation of the variance of the estimator for a fixed cost are the following:
\begin{itemize}
\item[(i)] the functions $f_k$ so that they do not depend on the starting point in $\partial B_k$; 
\item[(ii)] the optimal values of the parameters $N$, $M$, $\{R_k\}_{k=1}^M $ and $\{P_k\}_{k=1}^{M+1}$ obtained in \citet{Lagnoux08} for the unidimensional case (i.e. $s=1$). More precisely, $N$ is related to $C_{M+1}$ and all the $R_k$'s are equal to a same value, say $R$ which depends on $N$ and $C_{M+1}$. Furthermore, in order to satisfy the tradeoff between a premature death of the algorithm ($R_k P_{k+1} \ll 1$) and a prohibitive cost ($R_k P_{k+1} \gg 1$), we need the condition $R_k P_{k+1}=1$. Then $M$ is fixed by the relation $R p^{1/(M+1)}=1$.
\end{itemize}
\end{prop}

As expected, the optimal choice consists in taking the thresholds $\partial B_k$ in such a way that $f_k$ is constant. This is consistent with the observations of Section~\ref{ssec:gen_setting}. Nevertheless, the difficulty lies in the evaluation of the importance function $f$ and so in the design of the thresholds. We will see in Section~\ref{sec:perturb} the impact of a non optimal choice on the variance and on the cost of the algorithm.

If for some $k$, the function $f_k$ is constant, given that $\gamma_k(f_k)=p$, we get  the following identity $f_k=p/\gamma_k(\mathbf{1})$. Moreover, it comes from the definition of $\mu_k$ and Equations~\eqref{op-carre} and~\eqref{eq:muk_gk} that $\Gamma_k(f_k)=g_{k-1}(1-g_{k-1})p^2/\gamma_k^2(\mathbf{1})$ and 
\[
\gamma_{k-1}\left(\Gamma_k(f_k)\right)=\frac{p^2}{\gamma_k(\mathbf{1})}\frac{\mu_{k-1}\left(g_{k-1}(1-g_{k-1})\right)}{\mu_{k-1}(g_{k-1})}.
\]
Besides by~\eqref{eq:dyngamma}, the function $f_{k-1}$ can be expressed as
$f_{k-1}= g_{k-1}p/\gamma_k(\mathbf{1})$. Moreover, if $f_{k-2}$ is also constant, then $P_{k-1}(g_{k-1})=\gamma_k(\mathbf{1})/\gamma_{k-2}(\mathbf{1})$, and after calculus
\begin{equation*}
\gamma_{k-2}\left(\Gamma_{k-1}(f_{k-1})\right)=p^2 \left[\frac{1}{\gamma_{k-1}(\mathbf{1})} \frac{\mu_{k-1}(G^2_{k-1})}{\mu^2_{k-1}(g_{k-1})}-\frac{1}{\gamma_{k-2}(\mathbf{1})}\right].
\end{equation*}
Finally, if all the functions $f_k$ are constant, then the functions $g_k$ are also constant: $g_k= \gamma_{k+1}(\mathbf{1})/\gamma_k(\mathbf{1})$ and as for the functions $\Gamma_{k+1}(f_{k+1})$:
\[\Gamma_{k+1}(f_{k+1})=\frac{p^2}{\gamma_{k}(\mathbf{1})}\left[\frac{1}{\gamma_{k+1}(\mathbf{1})}-\frac{1}{\gamma_{k}(\mathbf{1})}\right].\]

\begin{rem}\label{rm:Mir}
These results justify the choices done in the algorithm proposed in \citet{MSM09}. The authors assume that 
\[
\underset{B\to \infty}{\lim} \frac 1B \log p_B^s=-\gamma(s) \quad \forall s\notin A;
\] 
where $p_B^s$ represents the probability to reach the target event $A$ starting from $s$, $B$ the rarity parameter and $\gamma$ is a decreasing function. The algorithm consists in taking:
\begin{itemize}
\item the replication numbers (except the last one) all equal to $R$;
\item the number of thresholds $n_B$ equals to $\left\lfloor B\gamma(s)/\log(R) \right\rfloor$;
\item the frontier $l_k$ of the intermediate threshold $L_k$ equals to
\[\left\{x\in D\quad /\quad \gamma(s)-\gamma(x)=\frac{k}{B}\log{R} \right\} \quad k=0\ldots n_B;\]
\item the last replication number equals to $R'=\left\lfloor e^{B\gamma(s)-n_B\log{R}} \right\rfloor$.
\end{itemize}
In other words, the authors equal all the replication numbers (excepted eventually the last one), take the number of thresholds equal to the optimal one in \citet{Lagnoux08}. Finally they fix all the thresholds in such a way that the decreasing rate $\gamma(s)$ is uniform over the thresholds and the probability to reach the target set $A$ starting from the $k$-th threshold depends on $k$ but not on the starting point of the frontier $l_k$.
\end{rem}

\section{Sensitivity analysis: deletion of a threshold}\label{sec:form_it_var}

Now, we study the sensitivity of $\Var(\widehat{p}_{M+1})$ with respect to the number of thresholds. We assume that the thresholds have the optimal shape: the functions $f_k$ are constant. It amounts to work in the unidimensional setting. Optimally, the thresholds are such that all the transition probabilities are equal, but $p$ being unknown this value cannot be computed. Moreover, in practice, the freedom of the choice of the thresholds can be limited by physical constraints. Then we study the need of an intermediate threshold and derive a procedure to detect whether we shall keep it or not.

\subsection{Iterative expressions of variance and cost}

The goal of this section is to compare the variance and the cost of the estimator obtained with $M$ thresholds with the ones obtained in the same setting but deleting the $k$-th threshold (thus in a simulation with $(M-1)$ thresholds). In that we view, we reallocate the replication numbers as following:
\begin{itemize}
\item for any $j=1, \ldots, k-2$, $R_j$ stays unchanged;
\item $R_{k-1}$ is replaced by $\lambda_{k-1} R_{k-1} R_k$;
\item and for any $j=k,\ldots, M-1$, $R_j$ is modified in $\lambda_{j} R_{j+1}$. 
\end{itemize} 

For instance, we can decide to keep all the $R_j$'s unchanged so the replication numbers are $R_1, \ldots, R_{k-1}$, $R_{k+1} \ldots, R_M$, or to report the replication number of the $k$-th threshold on the $k-1$-th's, the replication numbers being $R_1, \ldots, R_{k-1}R_k, R_{k+1}, \ldots, R_M$. 

\begin{prop}\label{prop:var_supp}
The variance of the estimator $\widehat{p}_{M+1}$ with $M$ thresholds is the sum of the variance of the estimator $\widehat{p}^{(-k)}_{M}$ obtained by running the algorithm  with the $k$-th threshold deleted (thus with $(M-1)$ intermediate thresholds) and the contribution of the $k$-th threshold: 
\begin{equation*}
\begin{split}
\Var\left(\widehat{p}_{M+1}\right)= \Var\left(\widehat{p}^{(-k)}_{M}\right) &+  \frac{1}{r_{k-1}}\left(1-\frac{1}{\Lambda_{k-1}R_k}\right) \gamma_{k-1}\left(\Gamma_k(f_k)\right)\\ 
&+ \sum_{j=k}^M \frac{1}{r_j}\left(1-\frac{1}{\Lambda_{j-1}}\right)  \gamma_j\left(\Gamma_{j+1}(f_{j+1})\right),
\end{split}
\end{equation*}
where $\Lambda_p=\prod_{j=k-1}^p \lambda_j$.

Similarly, the cost $C_{M+1}$ given in~\eqref{eq:cout} is the sum of the cost $C_M^{(-k)}$, computed with $M-1$ intermediate thresholds, and the contribution of the $k$-th threshold:
\[C_{M+1}=C_M^{(-k)} + r_{k-1}\left[ \gamma_{k-1}(c_{k-1})-R_k\Lambda_{k-1}\gamma_{k-1}(\tilde c_{k-1})  \right] + r_k \gamma_k(c_k) +\sum_{j=k+1}^M r_j \gamma_j(c_j) \left(1-\Lambda_{j-1}\right),\]
where $\tilde c_{k-1}$ stands for the cost of a particle going from the $(k-1)$-th threshold to the $k$-th in an algorithm with $(M-1)$ levels.
\end{prop}

The free parameters of the new algorithm with $M-1$ intermediate thresholds are $\{\Lambda_{j-1}\}_{j=k}^{M}$ that can be chosen by keeping the cost constant: it is sufficient to take
\begin{equation}\label{eq:cond_cout}
\Lambda_{k-1}=\frac{\gamma_{k-1}(c_{k-1})+R_k \gamma_{k}(c_{k})}{R_k\gamma_{k-1}(\tilde c_{k-1})} \quad \text{ and } \quad  \Lambda_{j-1}=1, \; j=k+1,\ldots, M.
\end{equation}
With these values, the variance $\Var\left(\widehat{p}_{M+1}\right)$ becomes
\[\Var\left(\widehat{p}^{(-k)}_{M}\right) +  \frac{1}{r_{k-1}}\left(1-\frac{1}{\Lambda_{k-1}R_k}\right) \gamma_{k-1}\left(\Gamma_k(f_k)\right) + \frac{1}{r_k}\left(1-\frac{1}{\Lambda_{k-1}}\right)  \gamma_k\left(\Gamma_{k+1}(f_{k+1})\right).\]

\subsection{Is the $k$-th threshold useful?}

Now, the goal is to study the need for an intermediate threshold and derive a procedure to detect whether we shall keep it or not. 
More precisely, the $k$-th threshold will be deleted if the variance of $\widehat{p}^{(-k)}_{M}$ is lower than the one of $\widehat{p}_{M+1}$, i.e. if the contribution 
of the $k$-th threshold is positive. In order to get a tractable procedure, we will assume the following: 
\begin{description}
\item[(A1)] All the thresholds have the optimal shape. Then we are lead to a unidimensional algorithm ($s=1$), so the measures $\gamma_j$ and the functions $f_j$, $g_j$ are constant;
\item[(A2)] The cost $\tilde{c}_{k-1}$ between $\partial B_{k-1}$ and $\partial B_{k+1}$ in the algorithm without the $k$-th threshold is given by 
$\tilde{c}_{k-1} =c_{k-1}+c_k$. Notice that each $c_k=c(g_k)$ is constant by (A1).
\end{description}
With these assumptions, we get
\begin{equation*}
\Lambda_{k-1}=\frac{a_k}{R_k}+g_{k-1}(1-a_k) \quad \text{and} \quad  \Lambda_j=1, \; j=k,\ldots, M-1,
\end{equation*}
where $a_k \coloneqq c_{k-1}/(c_{k-1}+c_k)$. Now, plugging these values into the variance, we get 
\begin{equation*}
\text{Var}(\widehat{p}_{M+1})=\text{Var}(\widehat{p}^{(-k)}_{M})+ \frac{p^2 Q(g_{k-1})}{r_{k}\gamma_{k+1}g_{k-1}\left[a_k+R_kg_{k-1}(1-a_k)\right]},
\end{equation*}
where 
\[
Q(x)=-x^2R_k(R_k\beta-1)(1-a_k)+x\left[R_k(R_k\beta-1)(1-a_k)-a_k(R_k-1)\right]+(R_k-1)\beta a_k\\
\]
and $\beta$ is defined by
\begin{equation}\label{def:beta}
\beta\coloneqq g_{k-1}g_k=\frac{\gamma_{k+1}}{\gamma_{k-1}}\in [0,g_{k-1}].
\end{equation}

Notice that $\beta=\mathbb{P}(\tau_{k+1}<\infty \mid\tau_{k-1} <\infty)$ quantifies the hardness for a particle to go from $\partial B_{k-1}$ to $\partial B_{k+1}$ and so
$\beta$ does not depend on the deleted $k$-th threshold.

The sign of $Q$ in the corrective term is the opposite of the one of the following polynomial
\[
R(x)= x^2-(1-\alpha)x-\alpha\beta,\qquad \trm{with} \qquad \alpha=\frac{a_k(R_k-1)}{R_k(1-a_k)(R_k\beta-1)},
\]
at $x=g_{k-1} \in ]0,1[$.
Its discriminant is $\Delta=(1-\alpha)^2+4\alpha\beta$.

In practice, we start by realising a pre-run in order to estimate the unknown parameters $\gamma_{k-1}$ and $\gamma_k$ and thus $g_{k-1}$ and $\beta$. Then the 
procedure is the following.

\begin{enumerate}
\item  If $R_k\beta=1$: $Q(g_{k-1})=a_k(R_k-1)\left(\beta-g_{k-1}\right)\leqp 0$ and it is recommended to preserve the $k$-th threshold.
\item  If $R_k\beta>1$: $\Delta$ is strictly positive and $R$ has two roots of opposite signs, $x_k^-<0<x_k^+<1$: 
\begin{enumerate}
\item when $0<g_{k-1}<x_k^+$, the polynomial $Q$ is positive and it is recommended to delete the $k$-th threshold;
\item when $x_k^+<g_{k-1}<1$, the polynomial $Q$ is negative and it is recommended to preserve the $k$-th threshold.
\end{enumerate}
\item If $R_k\beta<1$ and $\Delta<0$: the polynomials $R$ and $Q$ are positive and it is recommended to delete the $k$-th threshold.
\item If $R_k\beta<1$ and $\Delta>0$: the polynomial $R$ has two roots $x_k^-<x_k^+$:
\begin{enumerate}
\item when $0<g_{k-1}<x_k^-$, the polynomials $R$ and $Q$ are positive and it is recommended to delete the $k$-th threshold;
\item if $x_k^-<1$, when $x_k^-<g_{k-1}<1$, the polynomials $R$ and $Q$ are negative and it is recommended to preserve the $k$-th threshold.
\end{enumerate}
\end{enumerate}

Now we focus on the simplified cost because analytical values may be obtained. 

\begin{prop}\label{prop:cost_simpl}
Considering the simplified cost, there is no interest to introduce a new threshold when  $\beta \geqp 1/9$. When $\beta < 1/9$, the optimal positioning minimising the variance for a fixed cost is given by
$g_{k-1}=(1-3 \beta)/2$. In that case, the optimal replication number is
\[R_k^\ast=\frac{2(1-5\beta)}{1-9 \beta^2}\left(1 + \sqrt{\frac{2(1-\beta)}{1-5\beta}}\right)\]
that decreases from $2(1+\sqrt{2})$ for $\beta=0$ to $3$ for $\beta=1/9$.
\end{prop}

\section{Sensitivity analysis: perturbation of a threshold}\label{sec:perturb}

In this section, we assume all the thresholds $\partial B_i$ optimal (i.e. $f_i$ constant) except $\partial B_k$. Thus  
\begin{align*}
f_k&= \frac{p}{\gamma_{k+1}(\mathbf{1})}g_k, \quad P_k(g_k)=\mu_k(g_k)\mu_{k-1}(g_{k-1}),
\end{align*}
and the variance is given by
\begin{align*}
\frac{\Var(\widehat{p}_{M+1})}{p^2}=\frac{1}{\gamma_k(\mathbf{1})}\left(\frac{1}{r_{k-1}}-\frac{1}{r_{k}}\right)\frac{\Var_{\mu_k}(f_k)}{\E^2_{\mu_k}(f_k)} + \sum_{i=0}^M \frac{1}{r_i  \gamma_i(\mathbf{1})}\frac{1 - \mu_i(g_i)}{\mu_i(g_i)}.
\end{align*}
With a pre-run of the algorithm, we estimate the values of $g_k(i)$ for $i=1,\ldots,r$ and thus $\mu_k(g_k)$.

Now we want to twist $\partial B_k$ in order to get closer to the optimal shape and to obtain a new function $f_k$ constant. Consequently, with this new threshold, all the functions $f_k$ 
become constant and thus also the new function $g_k$, as explained in Section \ref{ssec:opt}. Introducing the new threshold 
$ \partial \widetilde B_k$ implies that $\gamma_k$, $P_k$, $P_{k+1}$, $g_k$ and $g_{k-1}$ are changed accordingly and we will use a $\widetilde{}$~symbol to denote the new terms. 

Furthermore, in order to guaranty a slight perturbation of threshold $k$, we assume naturally that
\[
B_{k+1} \subset \widetilde B_k \subset  B_{k-1},
\]
which implies that $\gamma_{k+1}(1)\leqp \widetilde \gamma_{k}(1)\leqp  \gamma_{k-1}(1)$.
We also introduce two operators $E_k$ and $E_{k+1}$ defined by
\[
\widetilde{P}_k=P_k E_k, \quad \widetilde{P}_{k+1}=E_{k+1}P_{k+1},
\]
and such that $\widetilde{P}_k \widetilde{P}_{k+1}=P_k P_{k+1}$. So defined, $E_k$ (respectively $E_ {k+1}$) is an operator acting on $\widetilde{\mathcal{F}}_k$ (resp. $\mathcal F_k$) valued in $\mathcal F_k$ (resp. $\widetilde{\mathcal{F}}_k$). We have
\[
\widetilde{g}_{k-1} = \widetilde{P}_k(1)=P_k(E_k 1) \quad \trm{and} \quad \widetilde{g}_k = \widetilde{P}_{k+1}(1)=E_{k+1}(g_k).
\]

Let us remark that $\widetilde{g}_{k-1}$ and  $\widetilde{g}_k$ are constant and linked by the identity
$\beta=\widetilde{g}_{k-1}\widetilde{g}_k$. If we choose $E_{k+1}(i,j)=\delta_{ij}/a_{i}$ with $a_i=Kg_k(i)$ for some constant $K$, then $E_{k+1}(g_k)=1/K$ so that $\widetilde g_k=1/K$.
Moreover, since $\widetilde g_k=\gamma_{k+1}(\mathbf{1})/\widetilde \gamma_{k}(\mathbf{1})$, we get
\begin{align}\label{eq:value_K}
K=\frac{\widetilde \gamma_{k}(\mathbf{1})}{\gamma_{k+1}(\mathbf{1})}.
\end{align}

Furthermore, taking $E_{k}(i,j)=Kg_k(j)\delta_{ij}$ leads to $E_kE_{k+1}=Id$ (and we recover $\widetilde P_k \widetilde P_{k+1}=P_k P_{k+1}$).
Finally,
\[
\widetilde P_{k+1}(i,j) = \frac{1}{K} \frac{P_{k+1}(i,j)}{g_k(i)} \quad \trm{ and } \quad \widetilde P_k(i,j) = K g_k(j)P_{k}(i,j).
\]
As a consequence, if $Kg_k(i)>1$, $\widetilde P_{k+1}(i,j)< P_{k+1}(i,j)$ for any $j$ and $\widetilde P_{k}(l,i)< P_{k}(l,i)$ for any $l$. It remains to determine the optimal value of $K$ that will be done by keeping the total cost of the algorithm constant which translates in
\[
\widetilde c_{k-1} \gamma_{k-1}(\mathbf 1)+R_{k} \widetilde c_k K \gamma_{k}(g_k)
= \gamma_{k-1}(c_{k-1})+R_{k} \gamma_{k}(c_k),
\]
leading to 
\begin{align*}
K=\frac{1}{\widetilde c_k} \left\{\frac{\mu_k(c_k)}{\mu_k(g_k)}+\frac{1}{R_k \beta} [\mu_{k-1}(c_{k-1})-\widetilde c_{k-1}]\right\}.
\end{align*}
Notice that fixing the value of $K$ amounts to defining the value of $\widetilde \gamma_k(\mathbf 1)$ by equation~\eqref{eq:value_K}. Remark that if the cost function $c$ is constant and equal to 1, then the optimal value of $K$ reduces to
$K=1/\mu_k(g_k)$.

\paragraph{Numerical application} Considering a two-dimensional Ornstein-Uhlenbeck process, we illustrate a way to deform the shape of a threshold in order to obtain an iso-probability levels. To this end, we simulate the stochastic process defined by
$$\begin{cases}
\diff X_t=-\Lambda X_t \diff t+\sigma \diff W_t,\; t> 0\\
X_0=x \in \R^2\\
\end{cases}$$
where $\Lambda=\textrm{diag}(\lambda_1,\lambda_2)$ with $\lambda_1>\lambda_2>0$, $\sigma>0$ and $W$ is a two dimensional standard Brownian motion.

\begin{figure}[h!]
\centering
\includegraphics[height=5cm]{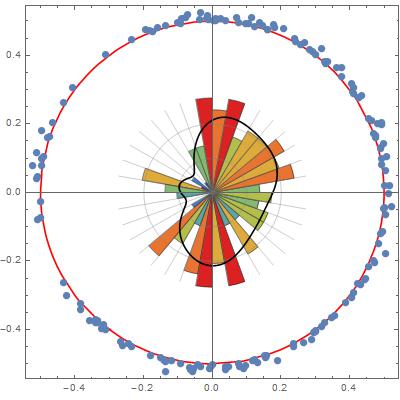}\hspace*{1cm}\includegraphics[height=5cm]{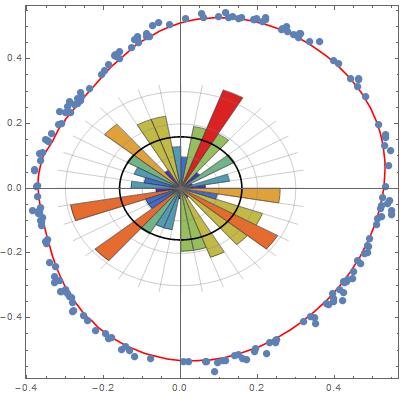}
\caption{\tiny The density of the occupancy measure at the first intermediate threshold and its estimation based on the von Mises kernel (black line). On the left, the threshold is the centered $0.5$ radius circle whereas on the right picture, the threshold is the conformal image of the circle.}\label{fig:dens1}
\end{figure}

We start the algorithm generating independently $N=300$ particles from $x=(0.05,0)$ and consider the $0.5$ radius circle as first intermediate threshold $\partial B_1$. In the sequel, we take $M=2$, $B_2=D(0,1)$ and $B=B_{M+1}=B_3=D(0,1.5)$. The parameters of the stochastic process are $\lambda_1=1, \lambda_2=0.2, \sigma= 0.3$ and its simulation is done via an Euler Scheme with a step of $0.01$ (we use the software Mathematica \cite{Mathematica15}). Firstly, we estimate the density of the occupancy measure of the process on $\partial B_1$\footnote{Since we work with continuous processes, the particles evolve until they reach $\partial B_1$ or the small disk $D(0,0.01)$ instead of the origin.}, with respect to its related Lebesgue measure. This estimation is based on the von Mises Kernel and as expected, this density (represented in Figure \ref{fig:dens1} left) is far from uniform. 

We previously noted that the efficiency of the splitting algorithm will be enhanced when the occupancy measures of the process on the intermediate thresholds are uniform. In our case, since $\lambda_1$ is greater than $\lambda_2$, we guess that the suited thresholds are ellipses. This intuition is confirmed by the left picture in Figure \ref{fig:dens1} and consistent with Theorem~(1.3) in \citet{ARG91} that establishes that, for any given $x\in \R^2$, $(Z_t)_{t\geqp 0}\coloneqq\left(\frac{\sqrt 2}{\sigma \sqrt{\log t}}X_t\right)_{t\geqp 0}$ admits the ellipse $\mathcal{E}=\{y=(y_1,y_2)\in \R^2;\, \lambda_1 y_1^2+ \lambda_2 y_2^2\leqp 2\}$ while $t$ goes to infinity.

\begin{figure}[h!]
\centering
\includegraphics[height=5cm]{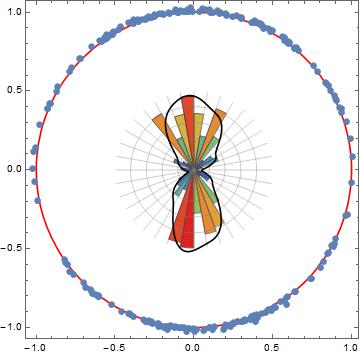}\hspace*{1cm}\includegraphics[height=5cm]{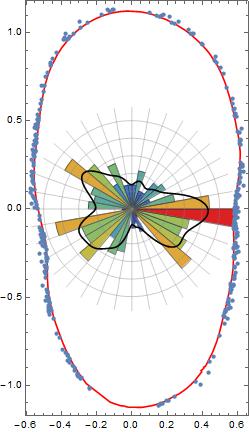}
\caption{\tiny Using a replication factor $R_1=2$ for the particles having reached the first deformed threshold, we make evolving these particles until they reach the next threshold or the inner $0.01$ radius circle. The empirical densities of the occupancy measure of the unit circle (left) and of the deformed threshold (right) and their respective estimations based on the von Mises kernel (black line) are represented. As mentioned before, we generate a new set of particles from the first deformed threshold instead of keeping the particles used to determine the conformal map. Thus in any rigor, we do not recover precisely the image measure; which explains the relative gap from the uniform distribution.}\label{fig:dens2}
\end{figure}

The goal is then to deform the first threshold. As the process lives in the plane, we can use a conformal map, $\varphi_1: B_1 \to \Omega_1$, in order to obtain a uniform occupancy distribution. Notice that the conformal maps are very convenient as planar transformations since they allow only local rotations and scales avoiding disturbing distortions; moreover, for common domains, $\partial \Omega_1=\varphi_1(\partial B_1)$. We follow the procedure described in \citet{WG10} to construct the conformal map (see also \ref{app:conformal} for more details). Once the conformal map $\varphi_1$ is computed (see the right picture in Figure \ref{fig:dens1}), we restart the algorithm using the threshold $\partial \Omega_1$ instead of $\partial B_1$. Then we start the next step which firstly estimates the density of the occupation measure on $\partial B_2$ after duplication of the particles in $\partial \Omega_1$ and secondly deform the shape of $\partial B_2$ as previously. Once $\partial \Omega_2$ is obtained, we restart the algorithm with this new threshold. Finally, the final step estimates the conditional probability to reach $\partial B_3$ for the particles in $\partial \Omega_2$. More precisely:

\begin{enumerate}
\item Each of the particles in $\partial \Omega_1$ is duplicated $R_1=2$ times and evolve independently from $\partial \Omega_1$ until $D(0,0.01)$ or $\partial B_{2}$ is reached. We determine the density of the occupancy measure of the process on $\partial B_{2}$ by using the von Mises kernel.
Then we find a conformal map $\varphi_{2}:B_{2}\to\Omega_{2}$ such that $\varphi_2(\partial B_2)=\partial \Omega_2$ and the image of the occupancy measure on $\partial B_{2}$ is the uniform measure on $\partial \Omega_2$ (See Figure~\ref{fig:dens2} for more details).

\item We perform independently a second set of particles from their same starting point at $\partial \Omega_1$ with the same size and stop them as soon as $D(0,0.01)$ or $\partial \Omega_2$ is reached.
\item Each of the particles in $\partial \Omega_2$ is duplicated $R_2=2$ times and evolve independently from $\partial \Omega_2$ until $D(0,0.01)$ or $\partial B_3$ is reached.  
\end{enumerate}

We emphasize that our intent is not to propose a new algorithm based on conformal mappings since we have not sufficient expertise to produce an efficient code. Working in 2D is already difficult, thus considering greater dimension becomes even more complex. Nonetheless, the harmonic functions or the quasi-conformal maps \cite{Ahlfors06,Heinonen06} are the natural generalization of the conformal transformations in higher dimensions. In our particular context, we start by estimating the density of occupation probability on a sphere what we can identify as a volume form. We can therefore attempt to determine a Riemannian metric $g$ such that associated Riemannian volume form is equal to the previous one. Then we can deform the metric $g$ into the uniform metric through a Ricci flow for instance; we get finally a new Riemannian variety homeomorphic to the sphere. For more details, see e.g. \cite{PLG14}.

\section{Conclusion and perspectives}

In this paper, we continue the multidimensional approach studied in (\citet{GHSZ98, Garvels00}) in order to obtain a new expression of the variance of the estimator analogous to that of the continuous case \citet{LLL09}. Then we derive the optimal parameters of the splitting algorithm. Furthermore, by introducing new operators, we obtain alternative expressions of the variance which are more tractable when we compare the variance of the estimators in an algorithm with $M$ thresholds and in an algorithm in which one of the threshold has been deleted. More precisely, we derive a procedure to detect whether we shall keep it or not. Finally, we investigate the sensitivity of the variance of the estimator with respect to a deviation of the threshold shape from the optimal one. We illustrate our theoretical results considering the planar Ornstein-Uhlenbeck process for which we propose a procedure based on conformal maps to twist the thresholds in order to get closer to the optimal shapes. 

A next natural research direction is probably the creation of a new algorithm that can decide the thresholds on the fly, for instance by using efficient algorithm for shape deformation. Such an algorithm has been proposed in \cite{Cerou-Guyader05} but only applies to 1D frameworks. When working in 2D or more, the problem is even more complex and challenging; see for instance \cite{DD07} for an approach based on the subsolutions of Hamilton-Jacobi-Bellman equations. Another way to investigate would consist in comparing, on a test example, the different existing algorithms dedicated to 2D or more (subset simulation \cite{AB01}, importance sampling based on design points \cite{APB99,DKD98} or adaptive pre-samples \cite{AB99}).

\section*{Acknowledgement} The authors are greatly indebted to the referees for their
fruitful and detailed suggestions or comments which permit us to greatly improve our paper.

\appendix

\renewcommand*{\thesection}{Appendix \Alph{section}} 
\section{The operator $\Gamma_k$ and its iterates}\label{app:op_Gamma}

The variance of the estimator of the target probability involves operators $\Gamma_k$ and their iterates defined in the following way. For $k=1,\ldots,M+1$, we introduce $\Gamma_k$ acting on $\mathcal F_k \times \mathcal F_k$ and valued in $\mathcal F_{k-1}$ by,
for $f$, $g\in \mathcal{F}_{k}$,
\begin{equation}\label{def:op_Gamma}
\Gamma_{k}(f,g)=P_{k}(fg)-P_{k}(f) P_{k}(g).
\end{equation}
To lighten notation, we denote $\Gamma_k(f)$ for $\Gamma_k(f,f)$. With the previous notation, one has 
\[\Gamma_1(f,g)=\gamma_1(fg)-\gamma_1(f)\gamma_1(g)\quad \trm{ and } \quad   \Gamma_{M+1}(f,g)=fg\times f_{M+1} (1-f_{M+1})=fg.\]

Firstly, we straightforwardly check that $\Gamma_k$ is bilinear and symmetric. Secondly, writing $\Gamma_k(f)(i)$ as a conditional variance, we get
\[\Gamma_k(f)(i)=P_k(f^2)(i)-P_k(f)(i)^2=\E^c\left(\left[f(X_k)\1_{\{\tau_{k}<\infty\}} - \E^c\left(f(X_k)\1_{\{\tau_{k}<\infty\}} \right) \right]^2\right)\]
where $\E^c$ is the expectation conditionally to the set $\left\{X_{k-1}=i;\; \tau_{k-1}<\infty\right\}$, we get
$\Gamma_k(f)\geqp 0$. Moreover, by~\eqref{def:G}, $\Gamma_k(f,\mathbf{1})=(1-g_{k-1})P_k(f), \quad \Gamma_k(\mathbf{1})=g_{k-1}(1-g_{k-1})$
and by~\eqref{eq:dyngamma},
\begin{equation}\label{op-carre}
\Gamma_k(f_k)=P_k(f_k^2)-f_{k-1}^2.
\end{equation}

%

Now let us iterate the construction of $\Gamma_k$. In that view, we introduce the multiplicative operator $\Gamma^{(0)}_k$ defined by 
\[\Gamma^{(0)}_k(f,g)=fg, \quad f,g \in \mathcal F_k,\]
such that $\Gamma_k(f,g)=P_k(\Gamma_k^{(0)}(f,g))-\Gamma_{k-1}^{(0)}(P_k(f),P_k(g))$.
This suggests to define for any $k$, the iterated operators $\Gamma_k^{(n)}$ in the following way:
\begin{equation}\label{op-carre2}
\Gamma_k^{(n+1)}(f,g)=P_{k-n}\left(\Gamma_k^{(n)}(f,g)\right) - \Gamma_{k-1}^{(n)}(P_k(f),P_k(g)), \quad k \geqp 1, \; 0 \leqp n \leqp k-1
\end{equation}
with the convention $P_1(f)=\gamma_1(f)$ precised above. The operator $\Gamma_k^{(n)}$ valued in $\mathcal F_{k-n}$ acts on $\mathcal F_k \times \mathcal F_k$.

We use the same simplified notation $\Gamma_k^{(n)}(f)$ to refer to $\Gamma_k^{(n)}(f,f)$. We introduce $P_{p,n}$ defined by
\begin{eqnarray}\label{defPnp}
P_{n,n}&=&Id, \quad n=0,\ldots, M,\\ \nonumber 
P_{p,n}&=&P_{p+1}\ldots P_n,\quad 0\leqp p \leqp n-1, \quad n=1,\ldots, M.
\end{eqnarray}
By induction on $n$, we easily get
\begin{equation}\label{op-carre3}
\Gamma_k^{(n)}(f) = \sum_{i=0}^n \binom{n}{i}(-1)^i P_{k-n,k-i}\left[\left(P_{k-i,k}(f)\right)^2\right].
\end{equation}
Since $P_{k-i,k}(f_k)=f_{k-i}$ and $f_0=p$, it comes in the particular case of $f=f_k$,
\begin{equation*}
\Gamma_k^{(n)}(f_k) = \sum_{i=0}^n \binom{n}{i}(-1)^i P_{k-n,k-i}\left(f^2_{k-i}\right).
\end{equation*}
Since for $f \in \mathcal F_k$, $\Gamma_k^{(n)}(f) \in \mathcal F_{k-n}$, we can compute $\gamma_{k-n}(\Gamma_k^{(n)}(f))$. From~\eqref{op-carre3} and the fact that $\gamma_{k_n} P_{k_n,k-i}=\gamma_{k-i}$, we get
\begin{equation*}
\gamma_{k-n}(\Gamma_k^{(n)}(f_k))=\sum_{i=0}^n \binom{n}{i}(-1)^i \gamma_{k-i}(f_{k-i}^2),
\end{equation*}
with $\gamma_0(f)=f$. $\Gamma_k^{(k)}(f_k)$ being a constant function in $\mathcal F_0$ equals $\gamma_0(\Gamma_k^{(k)}(f_k))$ and then the previous identity leads to
\[\Gamma_k^{(k)}(f_k)=\sum_{i=0}^k \binom{k}{i}(-1)^i \gamma_{k-i}(f_{k-i}^2).\]

The classical inversion formula which states the equivalence between the two following identities
\[u_k=\sum_{j=0}^k (-1)^j \binom{k}{j} v_j \quad  \textrm{and} \quad v_k=\sum_{j=0}^k (-1)^j \binom{k}{j} u_j\]
yields that
\begin{equation*}
\gamma_k(f_k^2) = \sum_{j=0}^k \binom{k}{j} \Gamma_j^{(j)}(f_j)
\end{equation*}
which means that $\gamma_k(f_k^2)$ can be written as the sum of terms involving the operators $\Gamma_k$ and their iterates. We get in particular the following identity
\begin{equation}\label{op-carre10}
\gamma_k\left(\Gamma_{k+1}(f_{k+1})\right)=\gamma_{k+1}(f_{k+1}^2)-\gamma_k(f_k^2)=\sum_{j=0}^k \binom{k}{j}\Gamma_{j+1}^{(j+1)} (f_{j+1}).
\end{equation} 
Actually this identity comes from a more general relation: first we make a change of parametrization in~\eqref{op-carre2} to get the following relation (valid for any function $f \in \mathcal{F}_{k+n}$),
\[P_{k+1}\left(\Gamma_{k+n}^{(n-1)}(f)\right)=\Gamma_{k+n}^{(n)}(f) + \Gamma_{k+n-1}^{(n-1)}(P_{k+n}(f)).\]
By a descendant induction on $p$, one gets for any $f \in \mathcal{F}_{k+1}$ and $0 \leqp p \leqp k$:
\begin{equation*}
P_{p,k}\left(\Gamma_{k+1}(f)\right)=\sum_{j=p}^{k} \binom{k-p}{j-p} \Gamma_{j+1}^{(j+1-p)}\left(P_{j+1,k+1}(f)\right).
\end{equation*}

If $f=f_{k+1}$, since $P_{j+1,k+1}(f_{k+1})=f_{j+1}$, we get 
\begin{equation}\label{op-carre9}
P_{p,k}\left(\Gamma_{k+1}(f_{k+1})\right)=\sum_{j=p}^{k} \binom{k-p}{j-p} \Gamma_{j+1}^{(j+1-p)}(f_{j+1}).
\end{equation}
It suffices to set $p=0$ to recover equation~\eqref{op-carre10} since $P_{0,k}(g)=\gamma_1 P_{1,k}(g)=\gamma_k(g)$ for any function $g \in \mathcal{F}_k$.
The use of~\eqref{op-carre} allows us to rewrite~\eqref{op-carre9} in the following way
\[P_{p,k+1}(f_{k+1}^2)-P_{p,k}(f_k^2)=\sum_{j=p}^{k} \binom{k-p}{j-p} \Gamma_{j+1}^{(j+1-p)}(f_{j+1})\]
and by a summation on $k$ from 0 to $p$, for all  $0 \leqp p \leqp k$, we get
\begin{align*}
P_{p,k+1}(f_{k+1}^2)-f_p^2&=\sum_{m=p}^k \sum_{j=p}^{m} \binom{m-p}{j-p} \Gamma_{j+1}^{(j+1-p)}(f_{j+1})
=\sum_{j=p}^{k} \binom{k-p+1}{k-j} \Gamma_{j+1}^{(j+1-p)}(f_{j+1}).
\end{align*}
When $p=0$, one gets
\[\gamma_{k+1}(f_{k+1}^2)-f_0^2=\sum_{j=0}^{k} \binom{k+1}{k-j} \Gamma_{j+1}^{(j+1)}(f_{j+1}),\]
which would have been also derived directly by a telescopic sum of~\eqref{op-carre10}.
The action of the measure $\gamma_p$ ($p \leqp k$) on~\eqref{op-carre9} leads to
\begin{equation*}
\gamma_{k}\left(\Gamma_{k+1}(f_{k+1})\right)=\sum_{j=p+1}^{k+1} \binom{k-p}{j-p-1} \gamma_p\left(\Gamma_j^{(j-p)}(f_j)\right).
\end{equation*}

This formula could be exploited to split the expression \eqref{eq:var14} of the variance $\Var(\widehat{p}_{M+1})$ in  two parts: 
\[
\sum_{i=0}^{l} \frac{1}{r_i } \gamma_i\left(\Gamma_{i+1}(f_{i+1})\right)+\gamma_l(D_{l+1,M})
\]
where $D_{l+1,M}$ is a quantity which depends only on the thresholds greater than $l$.

\section{Proofs}\label{app:proof}

\tbf{Proof of Proposition \ref{prop:var}} By the previous notation, Equations~\eqref{eq:est} and~\eqref{eq:2_somme}, the variance of the estimator can be written as
\begin{equation*}
\Var(\widehat{p}_{M+1})= \frac{1}{r_M^2}\Var(Z_{M+1})=\frac{1}{r_M^2}\sum_{i,j}\Cov\left(Y_{M+1}^{i},Y_{M+1}^{j}\right).
\end{equation*}
To compute the covariances in the right hand side of the previous equation, we use the classical formula
\begin{equation}\label{eq:decomp_cov}
\Cov(Y,Z)=\Cov(\mathbb{E}[Y|\mathcal{F}],\mathbb{E}[Z|\mathcal{F}]) +
\mathbb{E}[\Cov(Y,Z|\mathcal{F})]
\end{equation}
where $Y$ and $Z$ are two r.v.s and $\mathcal{F}$ a $\sigma$-algebra and
$\Cov(Y,Z|\mathcal{F})\coloneqq \mathbb{E}[YZ|\mathcal{F}]-\mathbb{E}[Y|\mathcal{F}]\mathbb{E}[Z|\mathcal{F}]$.

In our case, conditioning with respect to the $\sigma$-algebra generated by $\Z_M$ leads to, for any $(i,j)$,
\[\Cov(Y_{M+1}^{i},Y_{M+1}^{j})=R_M^2 f_{M}(i)f_{M}(j)\Cov(Z_{Mi}, Z_{Mj}) + R_M f_{M}(i)(1-f_{M}(i))\E[Z_{Mi}] \delta_{ij}.\]
The last term in the right hand side cancels for $i \neq j$ since conditionally to $\Z_M$ the  $r$ variables $Y_{M+1}^{i}$, $i=1,\ldots,r$ are mutually independent. Finally, introducing the covariance matrix $\Sigma_{n}(i,j)=\Cov(Z_{ni},Z_{nj})$  and using~\eqref{eq:moypar}, we derive the following expression
\[\Var(Z_{M+1})=R_M^2\|f_M\|^2_{\Sigma_M} + r_M \gamma_M(f_M(1-f_M)),\]
where $\|\cdot\|_{\Sigma_M}$ is the norm associated to the scalar product $\langle , \rangle_{\Sigma_M}$ defined by 
\[\langle f, g\rangle_{\Sigma_M}=\sum_{ij}f(i)g(j)\Sigma_M(i,j);\]
where $f$ and $g$ are two functions defined on $\{1,\dots,s\}$.

To compute the scalar product $\langle f, g\rangle_{\Sigma_M}$, we derive by induction the matrix $\Sigma_M$ and more generally the matrices $\Sigma_n$. The initial term $\Sigma_1$ is given by~\eqref{eq:decomp_cov} and can be rewritten as 
\[\Sigma_1(i,j)=\begin{cases}-N \gamma_1(i)\gamma_1(j), & i\neq j\\ N \gamma_1(i)(1-\gamma_1(i)),&i=j\end{cases};\]
and one gets
$\langle f, g\rangle_{\Sigma_1}=N\left(\gamma_1(fg)-\gamma_1(f) \gamma_1(g)\right)$. 
By Equation~\eqref{eq:somme}, we get $\Sigma_n(l,k)=\sum_{i,j} \text{Cov}(Y_{nl}^i, Y_{nk}^j)$ and conditioning by $\Z_{n-1}$, we have for $i=j$ to consider the two terms of the right hand side of~\eqref{eq:decomp_cov}; while for $i \neq j$, the last term cancels by conditional independence. The moment generating function of the random vector $\Y_n^i$, conditionally to $\Z_{(n-1)i}$ is given by
\[\varphi(t_1,\ldots, t_r)=\left[(1-g_{n-1}(i)) + \sum_{j=1}^r P_n(i,j) e^{t_j}\right]^{R_{n-1}Z_{(n-1)i}}.\]
By derivation of $\varphi$ (or using the multinomial distribution), we get directly that on one hand 
\[\E[Y_{nl}^i|\Z_{n-1}]=R_{n-1}P_n(i,l)Z_{(n-1)i}\]
and on the other hand
\[\text{Cov}(Y_{nl}^i, Y_{nk}^i|\Z_{n-1})=\begin{cases}
R_{n-1}P_n(i,l)(1-P_n(i,l))Z_{(n-1)i} & \quad k=l\\
-R_{n-1}P_n(i,l)P_n(i,k)Z_{(n-1)i} & \quad k \neq l.
\end{cases}\]
By Equation~\eqref{eq:moypar} and $\E[Z_{(n-1)i}]=r_{n-2}\gamma_{n-1}(i)$, we have
\[\text{Cov}(Y_{nl}^i, Y_{nk}^i)=\begin{cases}
R_{n-1}^2 P_n(i,l)^2 \Sigma_{n-1}(i,i) + r_{n-1} \gamma_{n-1}(i)P_n(i,l)(1-P_n(i,l))& \quad k=l\\
R_{n-1}^2 P_n(i,l)P_n(i,k) \Sigma_{n-1}(i,i) -r_{n-1} \gamma_{n-1}(i)P_n(i,l)P_n(i,k) & \quad k \neq l
\end{cases}\]
and for $i \neq j$,
$\text{Cov}(Y_{nl}^i, Y_{nk}^j)=R_{n-1}^2 P_n(i,l)P_n(i,k) \Sigma_{n-1}(i,j)$, that leads to the expression of $\Sigma_n(l,k)$ after a summation on $i$ and $j$.
Now 
\begin{align*}
\langle f, g\rangle_{\Sigma_n} & = \sum_{(N R_1\ldots R_{n-1}) \gamma_n(j).k,l}f(k)g(l)\Sigma_n(l,k)\\
&=R_{n-1}^2 \sum_{i,j,k,l} f(k)P_n(i,k) g(l)P_n(j,k) \Sigma_{n-1}(i,j) + A\\
&=R_{n-1}^2 \langle P_n(f), P_n(g)\rangle_{\Sigma_{n-1}} + A
\end{align*}
where
\begin{align*}
A&=r_{n-1}\sum_{i,l}\gamma_{n-1}(i) f(l)g(l)P_n(i,l)-r_{n-1}\sum_{i,l,k}\gamma_{n-1}(i)f(k)g(l)P_n(i,l)P_n(i,k)\\
&=r_{n-1}\gamma_{n-1} [P_n(fg)-P_n(f)P_n(g)]=r_{n-1} \gamma_{n-1}(\Gamma_n(f,g)).
\end{align*}
We are lead to the following induction relation
\begin{equation*}
\langle f, g\rangle_{\Sigma_n}=R_{n-1}^2\langle P_n(f), P_n(g)\rangle_{\Sigma_{n-1}} + r_{n-1} \gamma_{n-1}(\Gamma_n(f,g))
\end{equation*}
that, applied to the function $f_n$, yields 
$\|f_n\|^2_{\Sigma_n}=R_{n-1}^2 \|f_{n-1}\|^2_{\Sigma_{n-1}} + r_{n-1} \gamma_{n-1}(\Gamma_n(f_n))$, from which we deduce
\begin{align*}
\frac{\Var(Z_{M+1})}{r_M^2}&=\frac{1}{N}\left[\gamma_1(f_1^2)-\gamma_1^2(f_1)\right] + \sum_{i=0}^{M-2} \frac{1}{r_{M-(i+1)}} \gamma_{M-(i+1)}\left(\Gamma_{M-i}(f_{M-i})\right) + \frac{1}{r_M} \gamma_M[f_M(1-f_M)]\\
&=\frac{1}{r_0}\left[\gamma_1(f_1^2)-\gamma_1^2(f_1)\right] + \sum_{i=1}^{M-1} \frac{1}{r_i} \gamma_{i}\left(\Gamma_{i+1}(f_{i+1})\right) + \frac{1}{r_M} \gamma_M[f_M(1-f_M)].\\
\end{align*}

With the convention and Equation~\eqref{eq:p}, we get
\begin{equation}\label{eq:var0}
\frac{\Var(\widehat{p}_{M+1})}{p^2}=\sum_{i=0}^{M} \frac{1}{r_i \gamma_i(\mathbf{1})} \frac{\mu_i\left(\Gamma_{i+1}(f_{i+1})\right)}{\mu_i^2(f_i)}.
\end{equation}
Proceeding with the classical notation, valid for any probability $\mu$,
\[\E_\mu(f)\coloneqq \mu(f), \quad \Var_\mu(f)\coloneqq \mu(f^2)-\mu^2(f),\]
and using relation~\eqref{eq:muk_gk}, $\gamma_{M+1}(\mathbf{1})=p$ and $\gamma_0(\mathbf{1})=1$,
one gets the desired result.\QED

\tbf{Proof of Proposition \ref{prop:cost}} The cost of the first step of the algorithm (particles issued from $0$) is  
$N c_0=r_{0} \gamma_{0}(c_{0})$
and the one of the $n$-th step (particles issued from $\partial B_{n-1}$) for $n=2,\ldots,M+1$ is
\[\sum_{i=1}^s R_{n-1} Z_{(n-1)i}c_{n-1}(i).\]
Finally, Formula~\eqref{eq:moypar} leads to a mean total cost given by~\eqref{eq:cout}
since by convention $\gamma_0(c_0)=c_0$.\QED

\tbf{Proof of Proposition \ref{prop:opt}} The variance of the estimator is given by
\beq
\frac{\Var(\widehat{p}_{M+1})}{p^2}=\sum_{k=1}^M \frac{1}{\gamma_k(\mathbf{1})}\left(\frac{1}{r_{k-1}}-\frac{1}{r_{k}}\right)\frac{\Var_{\mu_k}(f_k)}{\E^2_{\mu_k}(f_k)} + \sum_{k=0}^M \frac{1 - \mu_k(g_k)}{r_k \gamma_{k+1}(\mathbf{1})}.
\eeq
The minimisation consists in a first step to cancel the terms (independent of the others)
\[\frac{\Var_{\mu_k}(f_k)}{\E^2_{\mu_k}(f_k)}\] 
which leads to take the functions $f_k$ constant on $B_k$ i.e. to require that the success probability from $\partial B_k^{(i)}$ does not depend on $i$. Then we are lead to the unidimensional setting and we fix $s=1$. $\gamma_k$ and $g_k$ are now real numbers between 0 and 1:
\begin{equation}\label{eq:opt_M}
\gamma_k\equiv\gamma_k(\mathbf{1})=\mathbb{P}(\tau_{k}<\infty) \quad \trm{and} \quad g_k\equiv\mathbb{P}(\tau_{k+1}<\infty\vert \tau_{k}<\infty).
\end{equation} 
In a second step, we minimise the other term of the  variance for a fixed cost. The variance and the cost can be rewritten in the following way
\[
\frac{\Var(\widehat{p}_{M+1})}{p^2}=\sum_{k=0}^M \frac{1 - \mu_k(g_k)}{r_k \gamma_{k+1}} \quad
\trm{and}  \quad C=\sum_{n=0}^{M} r_{n} \gamma_{n}c_{n}=N c_0+\sum_{n=1}^{M} r_{n} \gamma_{n}c_{n}.\]
From~\eqref{eq:opt_M}, we are lead to the optimisation problem with $s=1$ of \citet{Lagnoux08}. \QED

\tbf{Proof Proposition \ref{prop:var_supp}}
To compute the variance of the estimator $\widehat{p}^{(-k)}_M$ in the new setting, i.e. without the $k$-th threshold, we use formula~\eqref{eq:var0}. In particular, the $(k-1)$-th first terms are unchanged, while as we need to transport the function $f_{k+1}$ from $\partial B_{k+1}$ on $\partial B_{k-1}$, the $k$-th term becomes 
\[\frac{1}{\lambda_{k-1} r_k} \gamma_{k-1}(\tilde{\Gamma}_k(f_{k+1}))=\frac{1}{\lambda_{k-1} r_k}\gamma_k(\Gamma_{k+1}(f_{k+1})) + \gamma_{k-1}(\Gamma_k(f_k))\]
where $\tilde{\Gamma}_k(f_{k+1})=P_{k-1,k+1}(f^2_{k+1})-[P_{k-1,k+1}(f_{k+1})]^2$. Finally, the last terms are not modified except the replication numbers.

Defining $\Lambda_p=\prod_{j=k-1}^p \lambda_j$, the variance $\Var\left(\widehat{p}^{(-k)}_{M}\right)$ of the new estimator can be expressed as
\begin{equation*}
\sum_{j=0}^{k-2} \frac{1}{r_j} \gamma_j \left(\Gamma_{j+1}(f_{j+1})\right) 
+ \frac{1}{\Lambda_{k-1}r_k} \gamma_{k-1} \left(\Gamma_{k}(f_{k})\right)+ \sum_{j=k}^{M} \frac{1}{\Lambda_{j-1}r_j} \gamma_j \left(\Gamma_{j+1}(f_{j+1})\right)
\end{equation*} 
which leads to the result. 

In our context, all the $c_k$'s are equal to $1$, so the value of $\Lambda_{k-1}$ given by~\eqref{eq:cond_cout} becomes 
\begin{equation*}
\Lambda_{k-1}=\frac{1}{R_k} + \frac{\gamma_k(\mathbf{1})}{\gamma_{k-1}(\mathbf{1})}=\frac{1}{R_k} + g_{k-1}.
\end{equation*}
The variance is now given by
\begin{equation}\label{eq:var_iter}
\Var\left(\widehat{p}_{M+1}\right)=\Var\left(\widehat{p}^{(-k)}_{M}\right) +  \frac{p^2}{r_{k-1}\gamma_{k}(\mathbf{1})}(1-g_{k-1}) + \frac{p^2}{r_{k-1}\gamma_{k+1}(\mathbf{1})} S(R_k),
\end{equation}
where
$S(R_k)=\frac{1}{R_k}(1-g_k) - \frac{(1-\beta)}{1+R_k g_{k-1}}$
whose minimum is achieved at 
\[R_k^\ast=\left(\frac{1-g_k}{1-g_{k-1}}\right) \left(1 + \sqrt{\frac{1-\beta}{g_{k-1}-\beta}} \;\right)
\quad \textrm{ and } \quad
S(R_k^\ast)=-(g_{k-1}-\beta) \left[\sqrt{\frac{1-\beta}{g_{k-1}-\beta}}-1\right]^2.\]




%
The corrective term in formula~\eqref{eq:var_iter} of the variance rewrites, up to a positive multiplicative coefficient, as 
\[4 g_{k-1}^2(1-\beta) - (g_{k-1}-\beta)(1+g_{k-1})^2.\]
The sign of the previous expression is the same of the polynomial $R(x)=x^2-x(1-3 \beta) + \beta$ at $x=g_{k-1}$. So
\begin{itemize}
\item when $\beta>1/9$, $R(g_{k-1})$ is strictly positive;
\item when $\beta=1/9$, $R(g_{k-1})=(g_{k-1}-1/3)^2$ is positive and cancels at $1/3$;
\item when $\beta<1/9$, $R(g_{k-1})$ is minimum at $g_{k-1}^*=(1-3\beta)/2$ and $R(g_{k-1}^*)=\frac{1}{4}(1-\beta)(9\beta-1)<0$. This minimum decreases with $\beta$ from $0$ (for $\beta=1/9$) to $-1/4$ (for $\beta=0$).
\end{itemize}
The result now becomes obvious. \QED

\section{Finding the conformal map of Section \ref{sec:perturb}}\label{app:conformal}

The goal is to determine a conformal map $\varphi$ from a disk $B$ to a domain $\Omega\coloneqq\varphi(B)$ such that the image of the occupancy measure $\varphi_{*} m$ on $\partial B$ is the uniform measure on $\partial \Omega$ and  $\varphi(\partial B)=\partial \Omega$. First, we restrict $B$ to the unit disk. Since, for any Borel set $E\in \partial \Omega$,
\begin{align*}
\varphi_* m(E) &= m(\varphi^{-1}(E)) \overset{\textrm{def}}{=}  \int_{\partial B} \ind_{\varphi^{-1}(E)} h(\xi)d\xi
= \int_{\partial \Omega} \ind_E(\omega) h(\varphi^{-1}(\omega))\frac{d\omega}{|\varphi'(\varphi^{-1}(\omega))|}
\end{align*}
and we want $\varphi_{*} m (E) = \frac{1}{|\partial \Omega|} \int_{\partial \Omega} \ind_E(\omega) d\omega$, 
the conformal map $\varphi$ has to satisfy
$|\varphi'(\xi)|=h(\xi) |\partial \Omega|, \quad \forall \xi\in \partial B.$

Taking $|\varphi'(\xi)|=h(\xi)$ induces $|\partial \Omega|=1$. Since $\varphi$ is a conformal map, $\varphi'$ is holomorphic on $B$ and not null and 
$\log |\varphi'|=\log h$ is thus harmonic on $B$. Then we follow the procedure described in \citet{WG10}.

\begin{enumerate}
\item Since we work on the unit disk, we solve the Dirichlet problem and find its harmonic conjugate function concomitantly using the Schwarz integral formula \cite[Chap VII, \S 2]{Remmert91} that allows one to recover a holomorphic function, up to an imaginary constant, from the boundary values of its real part:
\[
\phi(z)=\int_0^{2\pi} \log h(e^{i\theta}) \frac{e^{i\theta}+z}{e^{i\theta}-z} \frac{d\theta}{2\pi}+ig(0),\; |z|<1.
\]
\item Now we consider $e^{\phi}$ which is holomorphic on $B$. Since $B$ is a simply connected set and taking the Cauchy integral, there exists a holomorphic function 
$\Phi$ on $B$ such that 
\[
\Phi(z)=\int_{[0,z]} e^{\phi(\omega)} d\omega,
\]
where $[0,z]$ is the segment that links 0 and $z$. Since $e^{\phi}$ never cancels, $\Phi$ is a conformal map. Thus we define $\varphi=\Phi$.
\end{enumerate}

In the case of a disk $B$ of radius $l$, we take $|\varphi'(\xi)|= h(\xi) 2\pi l$ instead of  
$|\varphi'(\xi)|=h(\xi)$ to get a boundary of length $2\pi l=|\partial \Omega|$.

\section*{Bibliography}

\bibliographystyle{elsarticle-harv}
\bibliography{bibliom}

%
%
%
\end{document}